\documentclass[11pt]{article}
\usepackage{cite}
\usepackage{mathrsfs}
\usepackage{amsfonts}
\usepackage{amsmath}
\usepackage{amsfonts,amssymb,color}
\usepackage{dsfont}
\usepackage{curves}
\usepackage{mathrsfs}
\usepackage{pifont}
\usepackage{amssymb}
\usepackage{hyperref}
\usepackage{latexsym,amsmath,amssymb,amsfonts,epsfig,graphicx,cite,psfrag}
\usepackage{eepic,color,colordvi,amscd}

\newtheorem{theorem}{Theorem}[section]

\newtheorem{lemma}{Lemma}[section]
\newtheorem{corollary}{Corollary}[section]

\newtheorem{claim}{Claim}[section]
\newtheorem{conjecture}{Conjecture}[section]

\newcommand{\qed}{\hfill\rule{0.5em}{0.809em}}

\def\emptyset{\mbox{{\rm \O}}}

\def\qed{\hfill \rule{4pt}{7pt}}

\def\pf{\noindent {\it Proof. }}

\begin{document}
	
\title{Perfect weight divisibility of (fork, co-fork$\cup K_1$)-free graphs}
	
	\author{Ran Chen\footnote{Email: 1918549795@qq.com }, \; Baogang Xu\footnote{Email: baogxu@njnu.edu.cn. Supported by 2024YFA1013902},  \; Miaoxia Zhuang\footnote{Corresponding author: 19mxzhuang@alumni.stu.edu.cn }\\\\
	\small Institute of Mathematics, School of Mathematical Sciences\\
	\small Nanjing Normal University, 1 Wenyuan Road,  Nanjing, 210023,  China}
\date{}

\maketitle

\begin{abstract}
	
	A {\em fork} is a graph obtained from $K_{1,3}$ (usually called {\em claw}) by subdividing an edge once, an {\em co-fork} is the complement graph of a fork, and a {\em co-cricket} is a union of $K_1$ and $K_4-e$. A graph is perfectly divisible if for each of its induced subgraphs $H$, $V (H)$ can be partitioned into $A$ and $B$ such that $H[A]$
	is perfect and $\omega(H[B]) < \omega(H)$. A graph is perfectly weight divisible if the same condition holds for any positive integer weights on its vertices, using the maximum weight of a clique instead of the clique number. Sivaraman proposed the conjecture that every fork-free graph is perfectly divisible in a paper by Karthick {\em et al.} [Electron. J. Comb. 29 (2022), P3.19.]. In the same paper, they proved that every ((\text{fork},\text{co-cricket}))-free graph is either claw-free or perfectly divisible. In this paper, we show that every (fork, co-fork$\cup K_1$)-free graph is perfectly weight divisible, and hence perfectly divisible. This improves some results of Karthick {\em et al.}.
	
	\begin{flushleft}
		{\em Key words and phrases:} fork-free graphs, chromatic number, clique number, perfect divisibility, perfect weight divisibility\\
			{\em 2020 Mathematics Subject Classification:}  05C15, 05C75\\
	\end{flushleft}
	
\end{abstract}

\section{Introduction}
All graphs considered in this paper are finite and simple. For $v\in V(G)$, let $N_G(v)$ be the set of vertices adjacent to $v$, and let $d_G(v)=|N_G(v)|$. For $X, Y\subseteq V(G)$, let $N_G(X)=\cup_{x\in X} N_G(x)\setminus X$ and $M_G(X)=V(G)\setminus (X\cup N_G(X))$. We use  $G[X]$ to denote the subgraph of $G$ induced by $X$, call $X$ a \textit{clique} if $G[X]$ is a complete graph, and call $X$ a {\em stable set} if $G[X]$ has no edges. The \textit{clique number} $\omega(G)$ of $G$ is the maximum size taken over all cliques of $G$. If it does not cause any confusion, we usually omit the subscript $G$ and simply write $N(v)$, $d(v)$, $M(v)$, $N(X)$ and $M(X)$. We say that $v$ is {\em complete} to $X$ if $X\subseteq N(v)$, and say that $v$ is {\em anticomplete} to $X$ if
$X\cap N(v)=\emptyset$. We say that $X$ is complete (resp. anticomplete) to $Y$ if each vertex of $X$ is complete
(resp. anticomplete) to $Y$.

If $1 \textless |X| \textless |V (G)|$ and every vertex of $V (G) \setminus X$ is either complete or anticomplete to $X$, then we
call $X$ a {\em homogeneous set}.  For $u,v\in V(G)$, we simply write $u\sim v$ if $uv\in E(G)$ and write $u\not\sim v$ otherwise. For $X,Y\subseteq V(G)$, write $N_X(Y)=N(Y)\cap X$ and
$N_X(v)=N(v)\cap X$.
For an induced subgraph $H$ of $G$, we abbreviate
$M(V(H))$ as $M(H)$.

We say that a graph $G$ contains a graph $H$ if $H$ is isomorphic to an induced subgraph of $G$. For a set ${\cal H}$ of graphs, we say that $G$ is ${\cal H}$-free if $G$ contains no $H\in {\cal H}$. If ${\cal H}=\{H_1,..., H_t\}$, we simply write $G$ is $(H_1,..., H_t)$-free instead.

Let $k$ be a positive integer. We say that $G$ is $k$-colorable if there is a mapping $\phi: V(G)\to\{1,2,...,k\}$ such that $\phi(u)\neq \phi(v)$ whenever $u\sim v$. The chromatic number $\chi(G)$ of $G$ is the smallest integer $k$ such that $G$ is $k$-colorable.

A family ${\cal G}$ of graphs is said to be \textit{$\chi$-bounded} if there is a function $f$ such that $\chi(G) \leq f (\omega(G))$ for every $G\in {\cal G}$, and if
such a function does exist for ${\cal G}$, then $f$ is called a \textit{binding function} of ${\cal G}$ \cite{bounded}. In addition, the class ${\cal G}$ is \textit{polynomially $\chi$-bounded} if it has a polynomial binding function. It has long been known that there are hereditary graph classes that are not $\chi$-bounded,
and it is also known that there are hereditary graph classes that are $\chi$-bounded but not
polynomially $\chi$-bounded \cite{0}. See \cite{CK2024,SR2019,survey} for more results and problems on this topic.

A \textit{fork} is the graph obtained from $K_{1,3}$ (usually called claw) by subdividing an edge once. The class of claw-free graphs is a subclass of fork-free graphs. It is known that there is no linear binding function even for a very special class of claw-free graphs. Kim \cite{Kim} showed that the Ramsey number $R(3,t)$ has order of magnitude $\Theta(\frac{t^2}{logt})$, and thus
there are claw-free graph $G$ with $\chi(G) \geq O(\frac{\omega^2(G)}{log\omega(G)})$. Chudnovsky and Seymour proved in \cite{claw-free} that every connected claw-free graph
$G$ with a stable set of size at least three satisfies $\chi(G) \leq 2\omega(G)$. Liu {\em et al.}\cite{LSWY}, confirmed a conjecture of Sivaraman (see \cite{KKS2022}), and  proved that $\chi(G) \leq 7\omega^2(G)$ for every fork-free graph $G$.

A \textit{hole} of $G$ is an induced cycle of length at least 4, and a $k$-hole is a
hole of length $k$. A $k$-hole is called an \textit{odd hole} if $k$ is odd, and is called an \textit{even hole}
otherwise. An \textit{antihole} is the complement of some hole. An {\em odd antihole}
is defined analogously. A graph $G$ is perfect if $\chi(H)=\omega(H)$ for each of its induced subgraphs $H$. In 2006, Chudnovsky {\em et al.} proved the {\em Strong Perfect Graph Theorem}.

\begin{theorem}\label{perfect}{\em \cite{CRST2006}}
	A graph $G$ is perfect if and only if $G$ is (odd hole, odd antihole)-free.
\end{theorem}

A \textit{perfect division} of a nonempty graph $G$ is a partition
of $V(G)$ into $A$ and $B$ such that $G[A]$ is perfect and
$\omega(G[B])<\omega(G)$.
A graph is \textit{perfectly divisible} if every nonempty induced
subgraph has a perfect division \cite{Hoang}.
Let $w:V(G)\to\mathbb{Z}_{>0}$ be a weight function.
For $X\subseteq V(G)$, write $w(X)=\sum_{x\in X}w(x)$.
For an induced subgraph $H$ of $G$, let $\omega_w(H)$ denote
the maximum weight of a clique in $H$, with
$\omega_w(\emptyset)=0$.
A partition $(A,B)$ of the vertex set of a nonempty induced
subgraph $H$ is a \textit{$w$-perfect division} if $H[A]$ is
perfect and $\omega_w(H[B])<\omega_w(H)$.
A graph $G$ is \textit{perfectly weight divisible}
\cite{2-divisible2019} if, for every positive integer weight
function $w$, every nonempty induced subgraph of $G$ has a
$w$-perfect division.
In particular, every perfectly weight divisible graph is
perfectly divisible, by taking $w\equiv 1$.

Sivaraman (see \cite{KKS2022}) proposed the following conjecture.
\begin{conjecture}\label{fork}{\em\cite{KKS2022}}
	The class of fork-free graphs is perfectly divisible.
\end{conjecture}

This conjecture is not even known to be true for claw-free graphs.

\begin{figure}[htbp]
	\begin{center}
			\includegraphics[width=8cm]{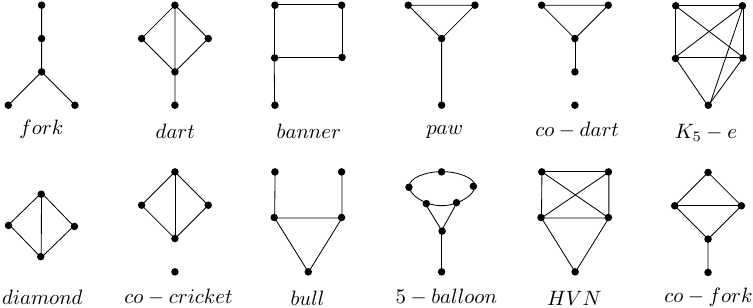}
		\end{center}
	\vskip -25pt
	\caption{Illustration of  fork and some forbidden configurations.}
	\label{fig-1}
\end{figure}

Let $G$ and $H$ be two vertex disjoint graphs. We use $G \cup H$ to denote a graph with vertex set $V (G) \cup V(H)$ and edge set $E(G) \cup E(H)$, and use $G + H$ to denote a graph with vertex set $V (G) \cup V (H)$ and edge set $E(G) \cup E(H) \cup \{xy | x \in V (G), y \in
V (H)\}$. An \textit{co-fork} is the complement graph of a fork.
A \textit{co-dart} is the union of $K_1$ and a paw. A \textit{bull} is a graph
consisting of a triangle with two disjoint pendant edges. A \textit{diamond} is the graph $K_1$ + $P_3$. A \textit{co-cricket} is the union of $K_1$
and a diamond. A \textit{balloon} is a graph obtained from a hole by identifying respectively two consecutive vertices with two leaves of $K_{1,3}$.
An \textit{$i$-balloon} is a balloon such that its hole has $i$ vertices. An $i$-balloon is called an \textit{odd balloon} if $i$ is odd. An {\em HVN} is a graph formed by adding a vertex that is adjacent to exactly two vertices of a $K_4$. See Fig. \ref{fig-1} for these configurations.

A graph $G$ is a \textit{minimal nonperfectly weight divisible graph}
if $G$ is not perfectly weight divisible but every proper induced
subgraph of $G$ is perfectly weight divisible.
Xu and Zhuang \cite{XZ2025} proved that every minimal nonperfectly
weight divisible fork-free graph is claw-free.
Consequently, the class of fork-free graphs is perfectly weight
divisible if and only if the class of claw-free graphs is perfectly
weight divisible.

Here, we list some results of perfect divisibility of fork-free graphs by further forbidding some graph $H$.
\begin{itemize}
	\item (Fork, HVN)-free graphs and ({\em fork}, $K_5-e$)-free graphs $G$ satisfy $\chi(G)\le \omega(G)+1$, and hence are perfectly divisible, see \cite{SR2019}.
	\item Karthick {\em et al.}  \cite{KKS2022} proved that Conjecture \ref{fork} is true on ({\em fork}, $F$)-free graphs when $F \in$ \{$P_6$, co-dart,bull\}. They also showed that every $G$ in the classes of ({\em fork}, $H$)-free graphs, when $H \in$ \{dart, banner, co-cricket\}, satisfies $\chi(G)\leq \omega^2(G)$.
	\item Wu and Xu\cite{WX2024} proved that Conjecture \ref{fork} is true on ({\em fork}, odd balloon)-free graphs. With minor modifications, the proof in \cite{WX2024} also shows that every ({\em fork}, odd balloon)-free graph is perfectly weight divisible.
	\item Xu and Zhuang \cite{XZ2025} proved that Conjecture \ref{fork} is true on ({\em fork}, $F$)-free graphs when $F \in$ \{$P_7$, $P_6\cup K_1$\}.
\end{itemize}

In this paper, we prove that every (fork, co-fork$\cup K_1$)-free graph is perfectly weight divisible and, consequently, perfectly divisible. Moreover, every such graph $G$ satisfies $\chi(G)\leq\binom{\omega(G)+1}{2}$. 
\begin{theorem}\label{co-fork}
	Every (fork, co-fork$\cup K_1$)-free graph is perfectly weight divisible.
\end{theorem}

\begin{corollary}\label{corollary-1}
	Let $G$ be a (fork, co-fork$\cup K_1$)-free graph. Then, $\chi(G)\leq\binom{\omega(G)+1}{2}$.
\end{corollary}

Notice that co-cricket is an induced subgraph of a co-fork$\cup K_1$. Co-cricket-free graphs are certainly (co-fork$\cup K_1$)-free. As a consequence of Theorem~\ref{co-fork}, (fork, co-cricket)-free graphs $G$ are perfectly divisible, and thus satisfy $\chi(G)\leq\binom{\omega(G)+1}{2}$.  These results improve some results of Karthick {\em et al.} in \cite{KKS2022} (they showed that if $G$ is (fork, co-cricket)-free, then $G$ is either claw-free or perfectly divisible, and $\chi(G)\leq\omega(G)^2$).

\section{Proof of Theorem~\ref{co-fork}}\label{main-23}
In this section, we prove the stronger statement that every
(fork, co-fork$\cup K_1$)-free graph is perfectly weight
divisible. Theorem~\ref{co-fork} then follows by
taking the all-ones weight function. The following lemmas
will be used.

\begin{lemma}\label{L-2}{\em\cite[Lemma~2.3]{WX2024}}
	Let $G$ be a minimal nonperfectly divisible {\em fork}-free graph. Then for each vertex $v \in V (G)$, $M(v)$ contains
	no odd antihole except $C_5$.
\end{lemma}

\begin{lemma}\label{L-3}{\em\cite[Corollary~3.1]{XZ2025}}
	Let $G$ be a minimal nonperfectly divisible fork-free graph, and $C = v_1 v_2 ... v_n v_1$ an odd hole contained in $G$. If there exist two adjacent vertices $u$ and $v$ in $V (G) \setminus V (C)$ such that $u$ is not anticomplete to $V (C)$ but $v$ is, then $N(u) \cap V(C) = \{v_i, v_{i+1}\}$, for some $i \in \{1, 2,...,n\}$.
\end{lemma}

\begin{lemma}\label{L-4}{\em\cite[Theorem~1.1]{XZ2025}}
Every minimal nonperfectly weight divisible fork-free graph
is claw-free.
\end{lemma}

Suppose to the contrary that $G$ is a minimal nonperfectly
weight divisible (fork, co-fork$\cup K_1$)-free graph.
Fix a positive integer weight function $w$ on $V(G)$ for
which $G$ has no $w$-perfect division. By minimality,
$G$ is connected, and by Lemma~\ref{L-4}, $G$ is claw-free.

Let $v_0\in V(G)$. If $G[M(v_0)]$ is perfect, then
$(M(v_0)\cup\{v_0\},N(v_0))$ is a $w$-perfect division:
indeed,
$\omega_w(G[N(v_0)])+w(v_0)\leq\omega_w(G)$.
Thus $G[M(v_0)]$ is not perfect. By Lemma~\ref{L-2}
and Theorem~\ref{perfect}, it contains an odd hole.
Choose an odd hole $C_0=v_1v_2\cdots v_nv_1$ in
$G[M(v_0)]$ with $n$ as small as possible.
Let $U=N(M(C_0))$. Clearly, $v_0\in M(C_0)$ and
$U\neq\emptyset$.
For $i\in \{1,2,\ldots, n\}$, let
\begin{eqnarray*}
	U_i&=&\{u\in U~|~N_{V(C_0)}(u)=\{v_i,v_{i+1}\}\}~\mbox{for}~1\leq i\leq n;\\	
	Z_i&=&\{x\in N(V(C_0))\setminus U~|~N_{V(C_0)}(x)=\{v_i, v_{i+1}\}\}; \\
	Z&=&\{x\in N(V(C_0))\setminus U~|V(C_0)\not\subseteq N(x)\};\\
	Z'&=&\{x\in N(V(C_0))\setminus U~|V(C_0)\subseteq N(x)\}.
\end{eqnarray*}
 All subscripts are taken modulo $n$.

It is easy to see that
$V(G)=V(C_0)\cup Z\cup Z'\cup U\cup M(C_0)$, and
$U=\bigcup_{i=1}^nU_i$ by Lemma~\ref{L-3}.
Since $G$ is claw-free, we have the following observation.
\begin{lemma}\label{lem-Z-NC}
For any $z\in Z$, $N(z)\cap V(C_0)\in \{\{v_i,v_{i+1}\}, \{v_i,v_{i+1},v_{i+2}\},\{v_i,v_{i+1},v_{j},v_{j+1}\}\}$ for some $i\in\{1,2,...,n\}$ and $j\in \{i+2,\ldots,i-2\}$.
\end{lemma}
\pf
Since $z$ is not complete to $V(C_0)$,
$G[N(z)\cap V(C_0)]$ is a disjoint union of paths.
If $v_i$ is an isolated vertex of $G[N(z)\cap V(C_0)]$, then
$\{v_i,v_{i-1},v_{i+1},z\}$ induces a claw, a contradiction. Moreover, since $G$ is claw-free, we have $\alpha(G[N(z)\cap V(C_0)])\le 2$.
Consequently, we have $G[N(z)\cap V(C_0)]$ is isomorphic to $P_2$, $P_3$, $P_4$,
or $2P_2$, which gives the assertion. \qed

\medskip
For $X_1,X_2\subseteq V(G)$, we use $[X_1,X_2]$ to denote the set of edges between $X_1$ and $X_2$. Next, we prove some properties of $U, U_i, Z, Z_i$ and $Z'$ in the  following (M1) to (M9).

\noindent{\bf (M1)} $Z=\bigcup_{i=1}^n Z_i$ and $[M(C_0),Z]=\emptyset$.

Let $z\in Z$ and $v\in M(C_0)$. By Lemma~\ref{lem-Z-NC},
either $z\in Z_i$ for some $i$, or
$N(z)\cap V(C_0)=\{v_i,v_{i+1},v_j,v_{j+1}\}$ for some
$i$ and $j\in\{i+3,\ldots,i-3\}$; otherwise,
$G[V(C_0)\cup\{v,z\}]$ contains a co-fork$\cup K_1$.
Suppose that the latter case occurs. Then $n\geq7$.
One of $v_{i+1}v_{i+2}\cdots v_jzv_{i+1}$ and
$v_iv_{i-1}\cdots v_{j+1}zv_i$ is an odd hole in
$G[M(v_0)]$ shorter than $C_0$, a contradiction.
Therefore, $z\in Z_i$ for some $i$, and hence
$Z=\bigcup_{i=1}^nZ_i$. This proves (M1).

\medskip

\noindent{\bf (M2)} $Z'$ is a clique anticomplete to $U\cup Z$.


First, we show that $Z'$ is anticomplete to $U\cup Z$.
Let $x\in U\cup Z$ and $z'\in Z'$. By symmetry,
suppose that $N(x)\cap V(C_0)=\{v_1,v_2\}$.
If $x\sim z'$, then $\{z',x,v_3,v_5\}$ induces a claw.
Thus $Z'$ is anticomplete to $U\cup Z$.

Since $U\neq\emptyset$, choose $u\in U_i$ for some $i$.
If $z_1',z_2'\in Z'$ are nonadjacent, then
$\{v_i,u,z_1',z_2'\}$ induces a claw.
Therefore, $Z'$ is a clique.

\medskip

\noindent{\bf (M3)} $U_i\cup Z_{i}$ is a clique.

Let $x,x'\in U_i\cup Z_i$. If $x\not\sim x'$, then $G[\{x,x',v_i,v_{i-1}\}]$ is a claw. Therefore, $x\sim x'$, and thus $U_i\cup Z_i$ is a clique.

\medskip

\noindent{\bf (M4)} $[Z_i,Z_k]=\emptyset$ for $k\notin\{i-2,i,i+2\}$, and $[Z_i,Z_{i+2}]$ is a matching.

Suppose $x\sim y$ for $x\in Z_i$, $y\in Z_k$ and $k\notin\{i+2,i,i-2\}$. If $k\notin\{i+1,i-1\}$, then $n\geq7$, and hence one of $xv_{i+1}v_{i+2}\cdots v_{k}yx$ and $xyv_{k+1}v_{k+2}\cdots v_nv_1v_2\cdots v_ix$ is an odd hole in $M(v_0)$ which has length less than $|V(C_0)|$, a contradiction to the choice of $C_0$. If $k\in\{i+1,i-1\}$, we may assume by symmetry that $k=i+1$, then $G[\{x,y,v_{i+1},v_{i+2},v_{i+3},v_0\}]$ is a co-fork$\cup K_1$.

Suppose that a vertex $x$ of $Z_i$ has two neighbors $y$ and $y'$ in $Z_{i+2}$. Since $y\sim y'$ by (M3), we have a co-fork$\cup K_1$ on $\{v_0, v_i, x, y, y', v_{i+3}\}$. Therefore, $[Z_i,Z_{i+2}]$ is a matching. This proves (M4).

\medskip

\noindent{\bf (M5)} If $z\in Z$ has a neighbor in $U_i$, then $z$ is complete to $\{v_i, v_{i+1}\}\cup U_i$.

Let $z\in Z$ and $u\in U_i$ such that $z\sim u$, and let $y\in N_{M(C_0)}(u)$. If $z$ has a nonneighbor in $\{v_i, v_{i+1}\}$, assume by symmetry that $z\not\sim v_i$, then $G[\{u, y, z, v_i\}]$ is a claw. If $U_i$ has a vertex $u'$ such that $u'\not\sim z$, then $z\sim v_{i-1}$ to avoid a claw on $\{u', v_i, v_{i-1}, z\}$, and  $z\sim v_{i+2}$ to avoid a claw on $\{u', v_{i+1},v_{i+2},z\}$. But now, $G[\{z,u,v_{i-1},v_{i+2}\}]$ is a claw. So, (M5) holds.

\medskip
By (M1), it also follows from (M5) that
$[U_i,Z_j]=\emptyset$ whenever $i\neq j$.
\medskip

\noindent{\bf (M6)} $N_{M(C_0)}(u)$ is a clique for any $u\in U$, and if $i\ne j$ and $u_i\sim u_j$ for $u_i\in U_i$ and $u_j\in U_j$, then $N_{M(C_0)}(u_i)=N_{M(C_0)}(u_j)$.

Without loss of generality, let $u\in U_i$ for some $i$. If $N_{M(C_0)}(u)$ has two nonadjacent vertices, say $x$ and $x'$, then $G[\{x,x',u,v_i\}]$ is a claw. Therefore, $N_{M(C_0)}(u)$ is a clique.

Suppose $i\ne j$, $u_i\in U_i$ and $u_j\in U_j$ such that $u_i\sim u_j$.  Since $i\ne j$, there exist two vertices $t_1,t_2\in \{v_i, v_{i+1}, v_j, v_{j+1}\}$ such that $t_1\in N(u_i)\setminus N(u_j)$ and $t_2\in N(u_j)\setminus N(u_i)$. If $N_{M(C_0)}(u_i)\ne N_{M(C_0)}(u_j)$, let by symmetry $y\in N_{M(C_0)}(u_i)\setminus N_{M(C_0)}(u_j)$, then $G[\{y, u_i, u_j, t_1\}]$ is a claw. Therefore, $N_{M(C_0)}(u_i)=N_{M(C_0)}(u_j)$. This proves (M6).

\medskip

 \noindent{\bf (M7)} For distinct vertices $x,y\in U_i$, $N_{M(C_0)}(x)\cap N_{M(C_0)}(y)=\emptyset$.

Suppose by symmetry that $x,x'\in U_1$ and $y\in N_{M(C_0)}(x)\cap N_{M(C_0)}(x')$. Then, $x\sim x'$ by (M3), and hence $G[\{x,x',y,v_2,v_3,v_5\}]$ is a co-fork$\cup K_1$. This proves (M7).

\medskip
If $z\in Z_i$ has two neighbors $x,y\in Z\setminus Z_i$,
then by (M4), we may assume that $x\in Z_{i-2}$ and
$y\in Z_{i+2}$.
But now, $G[\{z,x,y,v_i\}]$ is a claw.
Therefore, the edges between distinct sets $Z_i$ form a matching.

Let $(Z^{(1)},Z^{(2)})$ be a partition of $Z$ such that
each edge of this matching has one end in each part.
Then $[Z^{(r)}\cap Z_i,Z^{(r)}\setminus Z_i]=\emptyset$
for all $i$ and $r\in\{1,2\}$.
By (M3) and (M5), for each $z\in Z^{(r)}$,
$N(z)\cap(V(C_0)\cup U\cup Z^{(r)})$ is a clique.
\medskip

\noindent{\bf (M8)}
$\omega_w(G[U\cup Z^{(2)}])<\omega_w(G)$.

We show that every nonempty clique $K$ of
$G[U\cup Z^{(2)}]$ can be extended to a larger clique
of $G$. If $K$ contains a vertex of $Z^{(2)}\cap Z_i$,
then the choice of $Z^{(2)}$ and (M5) imply that
$K\subseteq U_i\cup Z_i$. Thus $v_i$ is complete to $K$.

Suppose that $K\subseteq U$. Let $\{i_1, i_2, \ldots, i_h\}\subseteq \{1, 2, \dots, n\}$ such that $K\cap U_{i_t}\ne\emptyset$ for $1\leq t\leq h$. If $h\ge 2$, then by (M6), $N_{M(C_0)}(u)=N_{M(C_0)}(v)$ for any two vertices $u,v\in K$. So, $M(C_0)$ has at least one vertex complete to $K$. If $h=1$, then $K\subseteq U_i$ for some $i$, then $v_i$ is complete to $K$. 
Since all weights are positive, (M8) follows.

\medskip

\noindent{\bf (M9)} $[U_i,U_{i+1}]=\emptyset$.

Suppose (M9) does not hold. Let, by symmetry, $x\sim y$ for some $x\in U_1$ and $y\in U_{2}$. If $n\geq7$, then $G[\{x, y, v_{2},v_{3},v_{4},v_{6}\}]$ is a co-fork$\cup K_1$. So, we suppose that $n=5$. If $M(C_0)$ is not a clique, then by (M6), there exists $v\in M(C_0)$ such that $v\not\sim x$ and $v\not\sim y$, which forces a co-fork$\cup K_1$ on $\{x,y,v_{2},v_{3},v_{4},v\}$. Therefore, we further suppose that $M(C_0)$ is a clique. By (M7), $U_1=\{x\}$ and $U_2=\{y\}$. 

Let $A=\{x,y,v_1,v_3,v_4,v_5\}\cup Z^{(1)}\cup M(C_0)$ and $B=\{v_2\}\cup (U\setminus\{x,y\})\cup Z^{(2)}\cup Z'$. We will show that $(A,B)$ is a $w$-perfect division of $G$.

The graph $G[A\setminus Z^{(1)}]$ is obtained by gluing the even hole $v_1xyv_3v_4v_5v_1$ and the clique $M(C_0)\cup\{x,y\}$ along the edge $xy$.
It is therefore perfect. Every vertex of $Z^{(1)}$ has a clique neighborhood in
$G[A]$, by the choice of $Z^{(1)}$ and the fact that
$Z$ is anticomplete to $M(C_0)$. Thus $G[A]$ is perfect.

Next, let $K$ be a clique of $G[B]$ with $w(K)=\omega_w(G[B])$. 
If $K\cap Z'\neq\emptyset$, then by (M2),
$K\subseteq Z'\cup\{v_2\}$, and $v_1$ is complete to $K$.
If $K\cap Z'=\emptyset$ and $v_2\in K$, then
$K\setminus\{v_2\}\subseteq Z_1\cup Z_2$.
Since $[Z_1,Z_2]=\emptyset$ by (M4), either $v_1$ or
$v_3$ is complete to $K$.
In both cases, $w(K)<\omega_w(G)$.
Otherwise, $K\subseteq U\cup Z^{(2)}$, and
$w(K)<\omega_w(G)$ by (M8).
Therefore, $(A,B)$ is a $w$-perfect division of $G$,
a contradiction. This proves (M9).

\medskip

Since $G$ is minimal nonperfectly weight divisible and $M(C_0)$ is a proper subset of $V(G)$, $G[M(C_0)]$ has a $w$-perfect division $(A_0,B_0)$.
Thus $G[A_0]$ is perfect and 
$\omega_w(G[B_0])<\omega_w(G[M(C_0)])$.
We will deduce a contradiction by extending $(A_0,B_0)$ to a $w$-perfect division of $G$.

\begin{claim}\label{X}
Let $X\subseteq U\cup B_0$ be a clique with
$w(X)=\omega_w(G)$.
Then there exist distinct indices $i_1,\ldots,i_h$,
where $h\geq1$, and $u_{i_j}\in U_{i_j}$ such that
\begin{itemize}
\item no two of $i_1,\ldots,i_h$ are consecutive modulo $n$;
\item
$N_{M(C_0)}(u_{i_1})=\cdots=N_{M(C_0)}(u_{i_h})$;
\item
$X=\{u_{i_1},\ldots,u_{i_h}\}\cup N_{M(C_0)}(u_{i_1})$.
\end{itemize}
\end{claim}

\pf By (M8) and the choice of $B_0$, we have that $X\cap U\ne\emptyset$ and $X\cap B_0\ne\emptyset$. By (M7), $|X\cap U_i|\leq1$ all $i$. By (M9), $U_i$ is anticomplete to $U_{i+1}$ for every $i$, so the indices of the sets $U_i$ meeting $X$
are pairwise nonconsecutive modulo $n$. Write $X\cap U=\{u_{i_1},\cdots,u_{i_h}\}$.
Since $X\cap U\ne\emptyset$, we have that $h\ge 1$. By (M6), $N_{M(C_0)}(u_{i_1})$ is a clique, and $N_{M(C_0)}(u_{i_1})=\cdots=N_{M(C_0)}(u_{i_h})$. Since $w(X)=\omega_w(G)$ and $X\cap B_0\ne\emptyset$, we have that $X\cap B_0=N_{M(C_0)}(u_{i_1})$. This proves Claim~\ref{X}. \qed

\medskip

For extending $(A_0, B_0)$ to a $w$-perfect division of $G$, we need to define a special subset $S$ of $U$. If $\omega_w(G[U\cup B_0])<\omega_w(G)$, then let $S=\emptyset$. Otherwise, let $X_1,X_2,\cdots,X_k$ be all the distinct
cliques contained in $U\cup B_0$ with
$w(X_j)=\omega_w(G)$. For each $j$, choose one vertex $s_j\in X_j\cap U$,
and let $S=\{s_1,\ldots,s_k\}$. It is certain that
\begin{equation}\label{S-1}
	\omega(G[(B_0\cup (U\setminus S)])<\omega(G).
\end{equation}

Next, we prove that
\begin{equation}\label{S-2}
	\mbox{$S$ is anticomplete to $A_0$.}
\end{equation}
Let $x\in S$, and let $X\subseteq U\cup B_0$ be a
clique containing $x$ with $w(X)=\omega_w(G)$. Let $\{i_1, i_2, \ldots, i_h\}$ be the set of integers satisfying Claim~\ref{X} with respect to $X$, and let $X=\{u_{i_1},u_{i_2},\cdots,u_{i_h}\}\cup N_{M(C_0)}(u_{i_1})$. It is certain that $x=u_{i_j}$ for some $j$. By the definition of $N_{M(C_0)}(u_{i_1})$, we have that $N_{M(C_0)}(\{u_{i_1},u_{i_2},\cdots,u_{i_h}\})=N_{M(C_0)}(u_{i_1})\subseteq B_0$. Therefore, $x$ is anticomplete to $A_0$. This proves (\ref{S-2}).

\medskip

Let $V_{odd}=\{v_1,v_3,v_5,\cdots,v_n\}$ and $V_{even}=\{v_2,v_4,\cdots,v_{n-1}\}$, and let
\begin{itemize}
	\item $A=A_0\cup V_{odd}\cup S\cup Z^{(1)}$;
		
	\item $B=B_0\cup (U\setminus S)\cup Z^{(2)}\cup V_{even}\cup Z'$.
\end{itemize}

It is easy to verify that $(A,B)$ is a partition of $V(G)$. The following two claims show that it is a $w$-perfect
division.

\begin{claim}\label{W-3}
	$\omega_w(G[B])<\omega_w(G)$.
\end{claim}
\pf
By (\ref{S-1}), $\omega_w(G[B_0\cup (U\setminus S)])<\omega_w(G)$. It suffices to prove that $\omega_w(G[B\setminus B_0])<\omega_w(G)$ since $B_0$ is anticomplete to $Z'\cup Z^{(2)}\cup V_{even}$.
Suppose that $X\subseteq B\setminus B_0$ is a clique
with $w(X)=\omega_w(G)$.
Since all weights are positive, $X$ is a maximal clique of $G$.

If $X\cap Z'\ne \emptyset$, then by (M2), $X\subseteq Z'\cup V_{even}$, and is complete to some vertex of $V(C_0)\setminus\{v_2,v_4,\cdots,v_{n-1}\}$ by the definition of $Z'$. Thus $X$ can be extended, a contradiction.
So, we suppose that $X\cap Z'=\emptyset$. 

By (M8), $X\cap V_{even}\ne \emptyset$. Let $v_2\in X$ by symmetry. Then 
$$X\subseteq \{v_2\}\cup (N(v_2)\cap B)=\{v_2\}\cup U_1\cup U_2\cup Z_1\cup Z_2.$$

By (M4), (M5) and (M9), there are no edges between
$U_{1}\cup Z_{1}$ and $U_{2}\cup Z_{2}$. Consequently, $X$ is complete to $v_1$ or $v_3$ since either $X\subseteq \{v_2\}\cup Z_1\cup U_1$ or $X\subseteq \{v_2\}\cup Z_2\cup U_2$. Therefore, $w(X)<\omega_w(G)$. This proves Claim~\ref{W-3}. \qed

\begin{claim}\label{P-3}
	$G[A]$ is perfect.
\end{claim}
\pf 
By (\ref{S-2}) and the definitions, $A_0$ is anticomplete
to $V_{odd}\cup S\cup Z^{(1)}$.
Moreover, every vertex of $Z^{(1)}$ has a clique
neighborhood in $G[A]$.
It therefore suffices to prove that $G[V_{odd}\cup S]$
is perfect. We first show that
\begin{equation}\label{e-6}
N_S(x)\subseteq U_i
\quad\hbox{for every }x\in S\cap U_i.
\end{equation}

Suppose $x\in S\cap U_c$ for some $1\leq c\leq n$. By the definition of $S$, $U\cup B_0$ has a maximum-weight clique, say $X$, of $G$ that contains $x$. By Claim~\ref{X}, we may let $X=\{u_{i_1},\cdots,u_{i_h}\}\cup N_{M(C_0)}(u_{i_1})$, where $u_{i_j}\in U_{i_j}$ for $1\leq j\leq h$, and $|i_j-i_k|\geq2$ (modulo $n$) for $1\leq j\ne k\leq h$, and $N_{M(C_0)}(u_{i_1})=\cdots =N_{M(C_0)}(u_{i_h})$. By symmetry, suppose $x=u_{i_1}$, i.e., $c=i_1$.

If there exists a $c'\not\in 
\{i_1, i_2, \ldots, i_h\}$ and a vertex $x'\in U_{c'}$ such that $x\sim x'$, then $c'\notin\{c-1, c+1\}$ by (M9), and $N_{M(C_0)}(x')=N_{M(C_0)}(x)=\ldots =N_{M(C_0)}(u_{i_h})$ by (M6). Hence, $x'$ is complete to $N_{M(C_0)}(u_{i_1})$. Since $X=\{u_{i_1},u_{i_2},\cdots,u_{i_h}\}\cup N_{M(C_0)}(u_{i_1})$ is a maximum-weight clique of $G$, $x'$ must have a nonneighbor in $\{u_{i_2},\ldots,u_{i_h}\}$, say $u_{i_2}$, by symmetry. But now, $G[\{x,x',u_{i_2},v_{i_1}\}]$ is a claw, a contradiction. Therefore,
$N_U(x)\subseteq U_c\cup\{u_{i_1},\ldots,u_{i_h}\}$.
Similarly, we have
\begin{equation}\label{e-7}
N_U(u_{i_\ell})\subseteq
U_{i_\ell}\cup\{u_{i_1},\ldots,u_{i_h}\}
\quad (1\leq\ell\leq h).
\end{equation}

If (\ref{e-6}) does not hold, then $x$ has a neighbor in $S\cap \{u_{i_1},u_{i_2},\cdots,u_{i_h}\}$, say $u_{i_{\ell}}$. This implies that $x,u_{i_{\ell}}\in S\cap X$. By the definition of $S$, one of $x$ and $u_{i_{\ell}}$ belongs to a maximum-weight clique of $G$ which is not $X$ but contained in $U\cup B_0$, contradicting (\ref{e-7}) and Claim~\ref{X}. Therefore, (\ref{e-6}) holds.

By (M3), we have that $U_i$ is a clique for $1\leq i\leq n$. Therefore, by (\ref{e-6}), for all $x\in S$, $N_{G[V_{odd}\cup S]}(x)$ is a clique. So, $G[V_{odd}\cup S]$ is perfect. This proves Claim~\ref{P-3}. \qed

By Claims~\ref{W-3} and \ref{P-3}, $(A,B)$ is a
$w$-perfect division of $G$, contrary to the choice of $w$.
Therefore, every (fork, co-fork$\cup K_1$)-free graph
is perfectly weight divisible, and hence perfectly divisible.
This completes the proof of
Theorem~\ref{co-fork}. \qed

\section*{Statements and Declarations}

\textbf{Conflict of interest} The authors declare that there are no known competing financial interests or personal relationships that could potentially influence the work presented in this paper.

\noindent\textbf{Data availability} Data sharing is not applicable to this paper since no datasets were generated or analyzed during the current study.

\end{document}